%% file: report_csc.tex
\documentclass[%
  onecolumn
   , colorlinks 
]{mpi2015-cscpreprint}

\usepackage{csquotes}
\usepackage[american]{babel}
\usepackage{amsmath, amssymb}
\usepackage{cleveref}
\usepackage{enumitem}

\input{mathmacros}

\input{config}

\begin{document}
\title{\thetitle}

\author[$\ast$]{\theauthori}
\affil[$\ast$]{\theaffiliationi{}.\authorcr%
  \email{\theemaili}, \orcid{\theorcidi}}

\author[$\dagger$]{\theauthorii}
\affil[$\dagger$]{\theaffiliationii{}.\authorcr%
  \email{\theemailii}, \orcid{\theorcidii}}

\author[$\diamond$]{\theauthoriii}
\affil[$\diamond$]{\theaffiliationiii{}.\authorcr%
  \email{\theemailiii}, \orcid{\theorcidiii}}

\author[$\ddagger$]{\theauthoriv}
\affil[$\ddagger$]{\theaffiliationiv{}.\authorcr%
  \email{\theemailiv}, \orcid{\theorcidiv}}

\shorttitle{\theshorttitle}
\shortauthor{\theshortauthor}
\shortdate{}


\msc{\themsc}

\abstract{\theabstract}

\novelty{A summary of perspectives from experts in matrix function computations is provided to help guide the field in its next steps.}

\maketitle

\begin{center}
	\itshape
	\thededication
\end{center}


\input{report_content}


\section*{Acknowledgments}%
\addcontentsline{toc}{section}{Acknowledgments}

\theacks


\addcontentsline{toc}{section}{References}
\bibliographystyle{abbrvurl}
\bibliography{references}

\end{document}

%% file: mathmacros.tex
\usepackage{bm, mathrsfs}
\usepackage[mathscr]{euscript}

\newcommand{\vb}{\bm{b}}

\newcommand{\vr}{\bm{r}}

\newcommand{\vu}{\bm{u}}
\newcommand{\vv}{\bm{v}}

\newcommand{\vx}{\bm{x}}



%% file: config.tex
\usepackage{todonotes}
\newcommand{\MF}{Max}
\usepackage{tcolorbox}
\tcbuselibrary{skins}
\tcbuselibrary{breakable}
\tcbset{shield externalize}
\newenvironment{mf}{%
	\list{}{\leftmargin0.2cm\rightmargin\leftmargin}\item\relax%
	\begin{tcolorbox}[breakable, %
		colback=Red!10!white, sharp corners, boxrule=0pt,frame hidden, %
		enhanced, borderline west={1mm}{-2mm}{Green}]%
		\small\sffamily{\color{Green}\textbf{\MF:}}}
	{\end{tcolorbox}\endlist}

\newcommand{\KL}{Kathryn}
\newenvironment{kl}{%
	\list{}{\leftmargin0.2cm\rightmargin\leftmargin}\item\relax%
	\begin{tcolorbox}[breakable, %
		colback=Red!10!white, sharp corners, boxrule=0pt,frame hidden, %
		enhanced, borderline west={1mm}{-2mm}{Violet}]%
		\small\sffamily{\color{Violet}\textbf{\KL:}}}
	{\end{tcolorbox}\endlist}

\def\showcomments{1}
\ifx\showcomments\undefinded
\usepackage{environ}
\NewEnviron{killcontents}{}

\fi

\newcommand{\theauthori}{Massimiliano Fasi}
\newcommand{\theorcidi}{0000-0002-6015-391X}
\newcommand{\theaffiliationi}{%
School of Computing,
University of Leeds,
Woodhouse Lane,
Leeds LS2 9JT,
United Kingdom}
\newcommand{\theemaili}{m.fasi@leeds.ac.uk}

\newcommand{\theauthorii}{Stéphane Gaudreault}
\newcommand{\theorcidii}{0000-0002-4475-0845}
\newcommand{\theaffiliationii}{%
Recherche en prévision numérique atmosphérique,
Environnement et Changement climatique Canada,
2121 Route Transcanadienne,
Dorval, Québec, H9P 1J3,
Canada}
\newcommand{\theemailii}{stephane.gaudreault@ec.gc.ca}

\newcommand{\theauthoriii}{Kathryn Lund}
\newcommand{\theorcidiii}{0000-0001-9851-6061}
\newcommand{\theaffiliationiii}{%
Max Planck Institute for Dynamics of Complex Technical Systems,
Sandtorstr.~1,
29106 Magdeburg,
Germany}
\newcommand{\theemailiii}{lund@mpi-magdeburg.mpg.de}

\newcommand{\theauthoriv}{Marcel Schweitzer}
\newcommand{\theorcidiv}{0000-0002-4937-2855}
\newcommand{\theaffiliationiv}{%
School of Mathematics and Natural Sciences,
Bergische Universität Wuppertal,
42097 Wuppertal,
Germany}
\newcommand{\theemailiv}{marcel@uni-wuppertal.de}

\newcommand{\thetitle}{Challenges in computing matrix functions}
\newcommand{\theshorttitle}{Challenges in $\lowercase{f}(A)\lowercase{\vb}$}
\newcommand{\theshortauthor}{M.~Fasi et al.}

\newcommand{\thededication}{This work is dedicated to the memory of Nicholas J.~Higham (1961--2024).}

\newcommand{\theabstract}{This manuscript summarizes the outcome of the focus groups at \textit{The f(A)bulous workshop on matrix functions and exponential integrators}, held at the \textit{Max Planck Institute for Dynamics of Complex Technical Systems} in Magdeburg, Germany, on 25--27 September 2023.  There were three focus groups in total, each with a different theme: knowledge transfer, high-performance and energy-aware computing, and benchmarking.  We collect insights, open issues, and perspectives from each focus group, as well as from general discussions throughout the workshop.  Our primary aim is to highlight ripe research directions and continue to build on the momentum from a lively meeting.}

\newcommand{\theacks}{The workshop itself was funded in part by DFG Project Number 529315380.  We thank all of the workshop participants (excluding the present authors) for their contributions to the substance of this manuscript:
Francesca Arrigo, Michele Benzi, Kai Bergermann, Philipp Birken, Liam Burke, Marco Caliari, Benjamin Carrel, Fabio Cassini, Ranjan Kumar Das, Vladimir Druskin, Andreas Frommer, Oswald Knoth, Patrick Kürschner, Thomas Mach, David Persson, Helmut Podhaisky, Michele Rinelli, Jonas Schulze, Roger Sidje, Igor Simunec, Martin Stoll, Mayya Tokman, Manuel Tsolakis, Paul Van Dooren, and Yannis Voet.}

\newcommand{\thekeywords}{exponential integrators, high-performance computing, matrix functions, matrix function times a vector, numerical linear algebra, research data management, workshop}

\newcommand{\themsc}{%
15A16, 
15-11, 
65F60, 
65Y05, 
65Y20, 
65-11  
}

%% file: report_content.tex
	The ``$f(A)\vb$'' community consists of researchers who develop, study, or use computational methods for computing the action of a matrix function on one or more vectors. These vectors can be considered either as a sequence of individual vectors or as a concatenation (\textit{block vector}).  A scalar function $f$ can be evaluated at a square matrix $A$ in a natural way that preserves many interesting properties of $f$.  Formally, $f(A)$ can be defined by means of the Jordan canonical form of $A$, the Taylor series expansion of $f$, its Cauchy integral representation, or Hermite interpolation~\cite[Chapter~1]{high:FM}.  If $f$ is analytic in a region that contains the spectrum of $A$, then these definitions are all equivalent.
	
	In numerical linear algebra, the expression ``computing matrix functions'' denotes two very different tasks:
	\begin{enumerate}[label=T\arabic*.,ref=T\arabic*]
		\item $f(A)$, i.e., the evaluation of the function $f$ at the $m \times m$ matrix $A$, which will produce an $m \times m$ matrix; and\label{it:fa-t1}
		\item $f(A)\vb$, i.e., the computation of the action of $f(A)$ on the $m \times n$ matrix $\vb$, where $n \ll m$, which will produce an $m \times n$ matrix.\label{it:fa-t2}
	\end{enumerate}
	In theory, any algorithm applicable to \ref{it:fa-t1} can be used for \ref{it:fa-t2} prior to a matrix--matrix product with $\vb$.  In many practical applications, however, the matrix $A$ is large and sparse, and the dense matrix $f(A)$ becomes impossible to store explicitly.  Some applications also require a linear combination of the action of multiple functions on different vectors, rendering this approach intractable.
	
	A prime example of such types of problems is exponential time integrators \cite{minchev2005review,hoos10}, a class of numerical methods for solving ordinary differential equations (ODEs) of the form
	\begin{align}\label{eq:ode}
		\dfrac{d }{d t} u(t) = F\bigl(u(t)\bigr), \quad u(t_0) = u_0.
	\end{align}
	Differential equations in this form appear in many areas of natural and social sciences.  In the majority of applications, the variable $u(t)$ represents an unknown dynamic quantity, $t$ is the independent variable, and $F$ describes the dynamics of the system.  Different types of exponential integrators can be derived by optimizing the coefficients of an ``ansatz'' of the form
	\begin{equation}\label{eq:ansatz}
		\varphi_0(A)\vv_0 + \varphi_1(A)\vv_1 + \varphi_2(A)\vv_2 + ... + \varphi_p(A)\vv_p,
	\end{equation}
	where the so-called $\varphi$-functions can be defined by the Taylor series
	\begin{equation*}
		\varphi_k(A) = \sum_{i=0}^\infty \frac{A^i}{(i+k)!}.
	\end{equation*}
	In other words, the ansatz~\cref{eq:ansatz} is just a linear combination of exponential-like functions evaluated at $A$ that act on a set of vectors.  Computationally, evaluating $\varphi_{i}(A) \vv_{i}$ is the most expensive step in exponential time integrators, and it has therefore been the subject of a considerable amount research.
	
	\setcounter{section}{-1}
	\section{The workshop}
	\label{sec:workshop}
	
	In September 2023, members of the $f(A)\vb$ community met at the Max Planck Institute for Dynamics of Complex Technical Systems in Magdeburg, Germany, to attend the first ``f(A)bulous workshop on matrix functions and exponential integrators'',\footnote{\url{https://indico3.mpi-magdeburg.mpg.de/event/30/}} organized by Kathryn Lund, Stéphane Gaudreault, and Marcel Schweitzer. The event featured traditional-style scientific talks covering recent algorithmic advances as well as applications, but a significant portion of the two and a half days was reserved for moderated discussion sessions.
	
	For these sessions, the workshop participants split into three focus groups (each led by one of the organizers of the workshop) in order to consider a broad challenge the community is facing, and they spent one afternoon looking at the issue and assembling potential solutions. For each group, a list of key questions was provided to foster and stimulate a lively---but also focused---discussion.  For the guiding slides with questions, see \cite{focus_group_slides}.
	
	The next day, each group reported the main discussion points to the other attendees, and the floor was then opened for further comments and questions.  Several participants (in particular, Thomas Mach and Yannis Voet) took thorough notes on the discussions and shared them with other participants via the workshop website.
	
	The three challenges that were selected for these focus groups were as follows:
	\begin{enumerate}
		\item knowledge transfer between the ``$f(A)\vb$ community'' and researchers from other areas (moderated by Marcel Schweitzer),
		\item high-performance and energy-aware computing (moderated by Stéphane Gaudreault), and
		\item benchmark problems and FAIR comparisons (moderated by Kathryn Lund).
	\end{enumerate}
	
	The next three sections of this report summarize the main conclusions of the three focus groups and the main points that were raised in the discussion that ensued.
	
	\section{Knowledge transfer}
	\label{sec:knowledge-transfer}
	
	The first group looked at ways to ensure that the algorithms developed within the community can reach the end users who need them.
	
	As discussed before, \cref{it:fa-t1,it:fa-t2} are fundamentally different from a computational perspective, and they generally require different techniques.  In particular, when only $f(A)\vb$ is sought, computing $f(A)$ is typically an unduly expensive and totally unnecessary step.
	
	Anecdotal evidence suggests that the distinction between the two problems is not as clear outside the $f(A)\vb$ community, and that researchers that wish to compute $f(A)\vb$ in an application domain often rely on the simple (but extremely inefficient) approach of multiplying $\vb$ by $f(A)$ after having computed the latter explicitly.
	
	There are several reasons behind this phenomenon. For one, there is no recent, authoritative source summarizing the existing methods for evaluating $f(A)\vb$.  The only comprehensive survey of the literature~\cite{frsi08} is now over 16 years old and does not reflect the breadth of methods currently available; more recent surveys and a thesis \cite{Gut13, GutKL20, Sch15} only deal with Krylov subspace methods and focus especially on limited-memory scenarios.  The situation is similar for exponential integrators, where the last comprehensive survey~\cite{hoos10} is 14 years old. This is in stark contrast to the literature on $f(A)$, which boasts a book~\cite{high:FM},\footnote{Incidentally, in the preface of~\cite{high:FM}, the author states ``The problem of computing a function of a matrix times a vector, $f(A)\vb$, is of growing importance, though as yet numerical methods are relatively undeveloped''. Due to this growing importance, a lot of developments have taken place since then, which are of course not covered by the last ``$f(A)\vb$ survey''~\cite{frsi08}, which was published in the same year as the book~\cite{high:FM}.} a survey paper dedicated to computational aspects~\cite{hial10}, and even a survey of existing software, which is periodically updated (and also includes software for $f(A)\vb$)~\cite{dehi16, hide14, hiho20}.
	
	An additional obstacle in knowledge transfer is the absence of $f(A)\vb$ in the standard academic curriculum.  In effect, most practitioners learn about the topic either during their graduate studies or through self-study, and their learning is hindered by the lack of a comprehensive review or textbook.
	
	For most functions $f$ of interest, it is relatively easy to find robust and efficient software toolboxes to evaluate $f(A)$.  These are available for most programming environments and work out-of-the-box for a large selection of test problems.  The situation is very different in the $f(A)\vb$ case, as the performance of the algorithm depends on a number of factors, and, with the software currently available, a certain level of experience is needed to find the right combination of parameters for a given computation.
	
	In most cases, the missing cornerstone is a reliable stopping criterion.  A stopping criterion requires a way of measuring or bounding the approximation error, but at present there is a knowledge gap in this area: available results typically only apply to normal matrices and are restricted to specific classes of matrices (e.g., Kronecker sums~\cite{besi16,besi15}) or specific functions (e.g., Cauchy--Stieltjes functions~\cite{fgs14a, frsc15, gukn13, gusc21, maro20}, Laplace transforms \cite{drus08,fkst23} or functions which can be related to an underlying ODE initial-value problem~\cite{bgh13,bokn20,bks23,bkt21}).  The error bounds available for more general cases may be very pessimistic, and they are typically not fit for use as stopping criteria in practice.  Furthermore, the derivation of error bounds quickly becomes interdisciplinary, as often function-specific, analytical results are necessary for deriving error expressions.  This is in stark contrast to, say, linear systems of the form $A \vx = \vb$, where a notion of residual $\vr := \vb - A \widetilde{\vx}$ for an approximation $\widetilde{\vx}$ is readily available from the data, can be cheaply approximated in Krylov subspace methods like GMRES, and is widely used as a reliable stopping criterion \cite{saad03}.
	
	The situation is further complicated by the lack of comprehensive comparisons of the performance of different algorithms, which is made arduous by the large number of implementation details and parameters that have to be selected.  We address this point further in Section~\ref{sec:benchmarking}.
	
	This focus group also identified a few timely challenges for the community, which represent, at the same time, opportunities for research advances.
	
	\begin{itemize}
		\item Current methods may not exploit the full breadth of available techniques, a case in point being randomized methods, which have just started to be considered~\cite{bugu23,ckn22,gusc23,pss23,pekr23}.
		\item Existing implementations are not able to fully leverage the variety of hardware on modern computers and supercomputers, which makes them potentially less attractive.  See also Section~\ref{sec:high-perf-energy}.
		\item Often, unlocking the full potential of methods might require ``intermingling'' the algorithm used for approximating $f(A)\vb$ with the surrounding ecosystem, taking into account specifics of the application at hand, instead of just treating it as a black-box that returns an approximate solution (see, e.g.,~\cite{chha23} for an example of a ``Krylov-aware'' approach in trace estimation, or~\cite{grt18} for an algorithm exploiting the intimate connection to exponential integrators when approximating linear combinations of $\varphi$-functions).
	\end{itemize}
	
	\paragraph{Next steps}
	
	Although a clear path to solve all the challenges affecting knowledge transfer is hard to pin down, some steps the community can take to make some progress on these issues are relatively easy to trace.
	
	First and foremost, there is a need for effective benchmarking standards.  All methods to be compared should be implemented from the same building blocks, and the same implementation choices should be applied consistently.  The performance of the methods should be measured in a uniform way, and this is not necessarily a simple task: estimating the execution time and memory usage of an algorithm is simple, but assessing its accuracy is not.  Accuracy is usually measured in terms of the forward error of the computed result, but since the exact solution is often not available, a reference solution computed using a different algorithm---potentially run in higher-than-working precision---is typically employed.  Attention should be paid to how the reference solution is computed, and how the error is estimated from it.
	
	Ideally, one could identify classes of problems where certain methods perform well, so that precise and easy-to-follow recommendations---based on the structure of $A$, the behavior of $f$, or a combination of both---can be made.
	
	In order for the comparison to be useful, it is also crucial that a representative set of benchmark problems be established.  Difficult questions that should be addressed regard:
	\begin{itemize}
		\item what information should be collected, in addition to the obvious $f$, $A$, and $\vb$, and
		\item what format should be used to store this information.
	\end{itemize}
	These points were discussed in more detail by another focus group---see Section~\ref{sec:benchmarking}.
	
	These efforts would put the community in a better position to summarize the literature and produce easy-to-follow guidance for all those interested in computing $f(A)\vb$ without delving into algorithmic and theoretical details: a general consensus is that drafting a modern survey of numerical methods for computing $f(A)\vb$ should be a priority.\footnote{Indeed, efforts in this direction are already underway.}  The lack of a comprehensive literature review for more specific problems, such as exponential integrators, is also a shared concern which should be addressed in coming years.
	
	Finally, a question that was raised is whether understanding better the sensitivity of the proposed algorithms, as well as offering ways of estimating the conditioning of a given problem, could help practitioners feel more confident about the use of a chosen algorithm for $f(A)\vb$---after all, this is one of the aspects that make the BLAS and LAPACK stand out.
	
	There does indeed exist a lot of work on estimating the condition number of the computation of $f(A)$ and $f(A)\vb$; see, e.g.,~\cite[Chapter~3]{high:FM} for a general overview of the topic. However, these condition number estimates are intimately related to the \emph{Fréchet derivative} $L_f(A,\cdot)$ of $f(A)$, an object that is typically several times more costly to compute than $f(A)$ itself. For the $f(A)$ problem, there exist several algorithms (each for a specific function $f$) that allow for the computation of $f(A)$ and its Fréchet derivative simultaneously and reuse certain computations in the process~\cite{AlHi09,AlHiRe13}. Building on this, \cite[Algorithm~7.4]{AlHi09} computes the matrix exponential $e^A$ together with a (quite reliable) condition number estimate at a cost of roughly 17 times that of computing $e^A$ alone.\footnote{The factor 17 can be reduced to 9 if a slightly reduced reliability is acceptable.}  Thus, even when using cleverly designed algorithms hand-tailored to a specific function, the overhead induced by the condition number estimator is quite substantial.
	
	Additionally, it is currently unclear how to extend such approaches to, e.g., Krylov subspace algorithms for $f(A)\vb$, as computing $f(A)\vb$ and $L_f(A,\cdot)$ has much less in common than computing $f(A)$ and $L_f(A,\cdot)$. Recently, there has been some progress in Krylov subspace algorithms for low-rank approximations of the Fréchet derivative~\cite{KaRe17,KKRS21,Kre19}, which might facilitate making a first step in this direction.

	\section{High-performance and energy-aware computing}
	\label{sec:high-perf-energy}
	
	The second focus group looked at the challenges surrounding the applications of $f(A)\vb$ in high-performance computing (HPC).  In many domains within the natural and social sciences and engineering disciplines, there is a need to compute $f(A)\vb$ where $A$ is extremely large and sparse.  Such large problems are typically solved using supercomputers, which are machines composed of many nodes with distributed memory, sometimes employing heterogeneous computing hardware. Each node is typically equipped with a number of CPUs and accelerators, such as GPUs, and to achieve peak performance, a routine must make the best use of all available resources.
	
	Aside from those for $A^{-1}\vb$, numerical methods to compute $f(A)\vb$ have seldom been used in large-scale parallel applications. One of the most active fields of research in this area is the solution of differential equations using exponential time integrators with Krylov and Leja point methods. Various factors impede the application of certain algorithms developed by the $f(A)\vb$ community, and we will provide a brief overview in this section.
	
	When the $A$ matrix is large and sparse, it is often impossible to store it explicitly in memory.  Fortunately, in many cases one can use ``matrix-free'' algorithms, which converge without the cost of forming or storing the matrix.  These are frequently used in HPC applications because they allow the solution of problems that would otherwise be intractable.  In the context of $f(A)\vb$, most matrix-free algorithms require only the action of the matrix (or an approximation to it) in the form of matrix--vector products.  For example, instead of storing the sparse Jacobian $J$ of a vector-valued function $F(\vu)$, its action on a vector can be approximated using the finite difference $J \vb \approx \left[ F(\vu + \epsilon \vb ) - F(\vu) \right] / \epsilon$, where $\epsilon$ is a small perturbation \cite{brown2008using}. Other matrix-free approaches, such as the complex-step approximation \cite{squire1998using} or automatic differentiation \cite{griewank2008evaluating}, are often also used.
	
	Requirements in terms of parallelism and memory storage considerably restrict the choice of possible algorithms.  The finite difference approximation of the Jacobian action illustrated above, for example, does not allow operations such as transposition, slicing, or pivoting without a prohibitive computational cost.  For the most part, it is not problematic to use a matrix--vector product routine instead of matrices for methods based on the Krylov subspace, the Taylor series, or Leja points. However, difficulties arise when information about the norm or the spectrum of $A$ is needed to 1) compute the parameters that make these methods efficient, or 2) determine the stopping criterion.  Without the matrix representation, it can be expensive to compute the operator norm or to estimate eigenvalues.  While the spectral radius $\rho(A)$ can be cheaply approximated using the power method on a single CPU, the communication cost in parallel implementations renders this idea inefficient in many HPC applications.  Further research will be necessary to develop numerical algorithms more suitable to this kind of problems.
	
	Another important consideration is the capability of an algorithm to scale and optimize energy efficiency as the amount of computing resources increases.  Clearly, existing implementations are not yet ready for the future exascale machines, i.e., supercomputers capable of performing at least $10^{18}$ binary64\footnote{Previously known as ``double precision''.} floating-point operations per second.  The problem of implementing Krylov subspace methods efficiently on a GPU has been considered from a theoretical point of view~\cite{fmyt16}, but studies focusing on high-performance implementations suggest that the use of GPUs is most beneficial when $A$ is dense~\cite{aeko18}, which is not often the case for $f(A)\vb$ problems, or when $A$ has a very specific sparsity pattern that can be mapped efficiently to GPU architectures~\cite{eios13}.  Therefore, software for this problem primarily targets CPUs and can only rely on GPUs in a limited number of cases or for a subset of the relevant operations.
	
	One obstacle is the fact that most implementations target binary64 accuracy, but binary64 arithmetic is not very efficient on GPUs.  When using the tensor cores on the latest NVIDIA H100 SXM5 GPUs~\cite[Table 1]{nvid22}, for example, the theoretical peak performance of the 19-bit TensorFloat-32 arithmetic is 494.7 trillion floating-point operations per second (TFLOPS), which improves for BFLOAT16 and binary16 arithmetic (989.4 TFLOPS) and breaks the PFLOPS barrier for the fp8 formats (1978.9 TFLOPS).  The peak performance of binary64 arithmetic is almost 30 times slower, with just 66.0 TFLOPS when tensor cores are used for matrix--matrix multiplications.  Using low-precision arithmetic is key to harnessing the full potential of GPUs, but only binary64 accuracy is typically sufficient for a range of applications.  In other areas of numerical linear algebra, this challenge has been addressed effectively by developing mixed-precision algorithms~\cite{aabc21, hima22}.  Efforts have been made recently in the context of exponential integrators to design such schemes \cite{balos2023leveraging}, but further research will be necessary.
	
	The issue is not solely with implementations: the algorithms themselves do not seem ready to address large-scale problems either. Many methods rely on matrix operations that do not scale well in a distributed environment. For example, full-basis orthogonalization, central to many Krylov subspace methods, necessitates numerous communication operations (i.e., message passing and synchronization) throughout the computation. This can result in unacceptable latencies, with most processes idly waiting for the slowest process to complete~\cite{BalCDetal14}.  There are a number of new developments on this front, including but not limited to low-synchronization orthogonalization \cite{BieLTetal22, CarLR21, CarLRetal22, FukKNetal20, FukNYetal14, OktC23, SwiLAetal21, YamTHetal20, Zou23a} and $s$-step methods \cite{Car20, CarGY22, YamCK22, YamTHetal20}, which can also be combined with one another.  Furthermore, a number of well established techniques are being rediscovered as communication-reducing, such as the natural short-term recurrences of Lanczos or batching vectors into tall-skinny matrices (block vectors) to take better advantage of BLAS Level 3 \cite{DreE19, DreE22}.  Sketching, randomization, and low precision can also be leveraged to reduce memory movement by shrinking the size of vectors to be stored and manipulated \cite{BalG22, BalG23, ckn22, gusc23, pss23, YamTKetal15}.  The performance of these techniques has been and is being thoroughly explored for linear systems solvers, but their transfer to matrix functions requires a better understanding of their backward stability, how they can be integrated into extended and rational Krylov subspace methods, as well as the conditioning of $f(A)\vb$ itself (cf.\ Section~\ref{sec:knowledge-transfer}).
	
	All these issues remain true for the computation of the full matrix $f(A)$ as well.
	
	A significant push in this direction could arise from an increased interest in exponential integrators. However, there are many factors that hinder their use in practical applications. Firstly, they are seldom featured in textbooks and are often absent from university curricula, resulting in many practitioners being unfamiliar with them. In addition, despite their better stability properties, they are generally more complicated to implement than explicit methods. Even when stability is important, practitioners tend to be more attracted by implicit or implicit-explicit schemes, because of the availability of techniques that deal with the stiffness of their particular problems. Furthermore, highly optimized libraries that implement algorithms for solving linear or nonlinear problems are readily available. This is not the case for exponential integrators, and the necessity to implement a parallel solver for the $\varphi$-functions often discourages the use of these methods in HPC applications. In recent years, there has been a renewed interest in exponential integrators, and this can largely be attributed to advances in numerical algorithms for the computation of $f(A)\vb$. While this is an encouraging trend, more work is needed to make these time integration schemes easier to use in applications.
	
	\paragraph{Next steps}
	
	Readying current $f(A)\vb$ work for exascale presents a number of significant challenges, but the community is well equipped to make some progress towards this goal.  It is clear that the focus should be on two distinct fronts, since not only the implementations, but also the algorithms, will require significant work in order to leverage the full computational power of next-generation supercomputers.
	
	In terms of rethinking existing algorithms, there is a clear need for reducing the number and frequency of communication operations. In particular, parallel inner products are a known communication bottleneck on distributed systems. Numerical algorithms with high arithmetic intensity should be favored over those requiring a high degree of data movements. Some work has already been done in this area (see, for example, block methods~\cite{FroLS17}, truncated orthogonalization~\cite{gusc23, kosk14, pss23} or restarts~\cite{aeeg08, bokn20, eier06, fgs14b}), but a completely different solution may be needed.
	
	In terms of implementations, a significant challenge is to leverage the untapped potential of GPUs, which, as they become faster and more prevalent, represent an increasingly large share of the overall performance of a supercomputer.
	
	Writing high-performance numerical linear algebra code that can target GPUs presents various difficulties. First and foremost, the variation in capabilities between different models of GPUs, especially those from different vendors, is dramatically larger than the variation between CPUs. As a consequence, there is no unified implementation of the BLAS and LAPACK for GPUs.  Vendors provide highly-optimized libraries for their own hardware, but these require very different frameworks, which means that porting an implementation from one GPU to another requires a significant human effort.  For example, NVIDIA provides cuBLAS,\footnote{\url{https://docs.nvidia.com/cuda/cublas/}} which is part of the CUDA Toolkit,\footnote{\url{https://developer.nvidia.com/cuda-toolkit/}} while AMD provides support through rocBLAS,\footnote{\url{https://rocm.docs.amd.com/projects/rocBLAS/}} which is part of the ROCm Platform.\footnote{\url{https://rocm.docs.amd.com}}
	
	Potential solutions, which include the C++ runtime API HIP,\footnote{\url{https://rocm.docs.amd.com/projects/HIP/}} also part of the ROCm Platform, the programming model SYCL\footnote{\url{https://www.khronos.org/sycl/}}~\cite{sycl20}, and libraries such as MAGMA\footnote{\url{https://icl.utk.edu/magma/}}~\cite{tdb10,tnld10}, are not yet mature enough to be used in production code.
	
	At present, the community should attempt to rewrite existing algorithms to ensure optimal performance on HPC architectures. They should seek to minimize communications and use low precision (binary32, binary16, or lower) for the bulk of the computation, switching to higher precision (typically binary64) only when strictly necessary.
	
	For research reproducibility, the $f(A)\vb$ community should adopt and promote open science best practices. This entails authors sharing the code and data that would allow to replicate the results presented in their publications.
	
	\section{Benchmarking}
	\label{sec:benchmarking}
	
	The last focus group discussed best practices for sharing code and data sets so that they can be easily reused, in accordance with FAIR guidelines.\footnote{\url{https://www.go-fair.org/}}  The main goal is to simplify two important steps of the algorithm development process:
	\begin{itemize}
		\item evaluating new implementations on established test problems, and
		\item comparing their performance with that of existing algorithms in the literature.
	\end{itemize}
	A welcome side effect, which comes at no additional cost, is the reproducibility of experimental results.  This idea promises to solve a number of problems that commonly arise when new algorithms are proposed in the literature.
	
	Unless the authors decide to compare their proposed new method with all state-of-the-art algorithms for the same problem, it is impossible for the reader to understand how the new method compares with existing alternatives.  A new implementation could easily perform worse than a much simpler and well established one, but the reader would have to spend a significant amount of effort to check whether this is the case, especially if the code used in the original publication is not available.
	
	The peer review process can help with this, but there are limitations.  Reviewers can recommend that new approaches be compared with the most relevant existing alternatives, but it is difficult to ensure that the comparison is fair, and most journals in numerical analysis and numerical linear algebra do not yet require submission of software or reproducibility of the experimental results.  Moreover, as a test set of representative $f(A)\vb$ problems is not currently available, it is difficult for a reviewer---and for the reader, later on---to make sure that the numerical experiments reported in a publication provide an impartial representation of the merits and drawbacks of new algorithms.  Not all methods are suitable for a given choice of $f$, $A$, and $\vb$, and having a battery of tests with clear classes of functions and matrices can help identify what types of problems a certain algorithm can deal with effectively.  This can help corroborate theoretical results, in addition to providing a quick and standardized way of comparing all relevant algorithms for a specific choice of $f$, $A$, and $\vb$.
	
	The metrics against which these algorithms should be compared are also not uniquely determined, and authors are free, within reason, to choose the ones that suit them best.  In some cases, the metric itself is poorly defined and can depend on a range of factors that are not within the control of who is performing the test.  A case in point is runtime, which is commonly used to assess the performance of different implementations on a same test set.  Runtime is very sensitive to the hardware configuration, as well as some low-level details of the software libraries being used, so that algorithm $\alpha_{1}$ can easily be faster than algorithm $\alpha_{2}$ on a machine and slower on another for the same test problem.
	
	When an algorithm cannot be implemented in the most efficient way possible, for example, because of limitations of existing hardware, an appealing alternative is to rely on the number of floating-point operations being performed.  This metric is only meaningful for large matrices, and it can be very inaccurate on modern hardware and especially in distributed-computing settings, as it focuses on arithmetic intensity when, in practice, the performance of most algorithms is bounded by the memory bandwidth.
	
	A final difficulty is represented by the lack of clear licensing for code and test problems alike.  This prevents reuse and, in many cases, hinders reproducibility of existing results.  For more information on research data management in mathematics, see a recent white paper by the Germany-based Mathematical Research Data Initiative (MaRDI) \cite{mardi23}.
	
	\paragraph{Next steps}
	
	It is a priority for the community to produce a set of representative test cases whereon new and old algorithms can be compared.  Ideally, one would want access to a remote facility that is capable of testing submitted implementations against known and unknown benchmark problems, providing overall scores for a number of metrics including accuracy, stability, and runtime performance.  Similar services exist for machine learning research~\cite{club21}, where the unknown problems are used to prevent authors from overfitting their models to the test set.
	
	The main difficulty to address in order to deliver this golden standard is to ensure a fair comparison among implementations written using different languages, as Julia, MATLAB/GNU Octave, and Python are all well established in this community, and being able to compare code in these different languages is likely to pose significant challenges in terms of software engineering.
	
	A more modest but attainable result would be the development of a curated reference collection of test problems.  Authors testing their code could then simply choose which parts of the collection to include, and they could justify their choices by pointing out which classes of problems are not suitable in their context.  Reviewers could equally rely on such a collection to ensure that authors are providing a fair picture of the merits of their algorithms.
	
	Building and maintaining this infrastructure would come with some logistical challenges.  A sufficient number of examples should be included, so that the collection represents the main applications in which $f(A)\vb$ appears.  As test examples can be quite large, the collection might require a hosting service with sufficient storage space and bandwidth, or standardized protocol to point to resources like Zenodo, from which data could be downloaded.  A possible remedy for the latter issue is to promote the use of so-called procedural examples, whereby a problem is specified mathematically and the matrix $A$ and vector $\vb$ can be generated with a desired size or other properties via a script.
	
	It is necessary to ensure that the test cases remain relevant, and that the collection grows and remains representative despite hardware and algorithm improvements that may make problems that are difficult today trivial in the near future.  The test cases will likely come from a number of researchers in various research domains, and will have to be collected and added to the collection by a number of volunteers.  A standard license---or set of licenses---should be adopted that ensure reproducibility and that fair credit is given to test problem creators and curators.
	
	Although it will not be possible for the community to enforce such a requirement, publishing and advertising a reasonable set of recommendations should be one of the priorities of the group working on this.
	
	\section{Conclusions}
	It is easy to feel overwhelmed looking at the long to-do lists we have outlined.  A change of perspective may lessen the anxiety: these are exciting opportunities, some of them even so-called ``low-hanging fruit'', and the impact of addressing them is huge, even for such a small field.  Matrix functions continue to surface in diverse applications, and many of the techniques developed for $f(A)\vb$ can cross-pollinate work in linear systems, matrix equations, Fréchet derivatives, and other problems we are not yet aware of.  Furthermore, the development of comprehensive surveys and language-agnostic benchmarking workflows for $f(A)\vb$ can set an example for other mathematical fields that are struggling to modernize and keep up with an ever-increasing publication load.  Our primary aim is that this manuscript builds on the momentum of a successful workshop and inspires new, meaningful projects in $f(A)\vb$ and beyond.